\documentclass{svproc}
\usepackage{url}
\usepackage{amsmath}
\usepackage{amsfonts}
\usepackage{subfigure}
\usepackage{multirow}
\usepackage{amssymb}
\usepackage{graphicx}
\usepackage[english]{babel}
\usepackage{fancyhdr}

\begin{document}
\mainmatter              
\title{Semi-analytical Solutions for Breakage and Aggregation-breakage Equations via Daftardar-Jafari Method}
\titlerunning{DJM for Breakage and Aggregation-breakage Equation}  
%
\author{Sanjiv Kumar Bariwal, Gourav Arora, and Rajesh Kumar }
\authorrunning{S.K.Bariwal et al.} 
%
%
\institute{Department of Mathematics, Birla Institute of Technology and Science, Pilani,\\ Pilani-333031, Rajasthan, India\\
\email{ p20190043@pilani.bits-pilani.ac.in}}

\maketitle              

\begin{abstract}
The semi-analytical method obtains the solution for linear/ nonlinear ODEs and PDEs in series form.
This article presents a  novel semi-analytical approach named Daftardar-Jafari method (DJM) to solve integro-partial differential equation such as breakage and nonlinear aggregation-breakage equations (ABE). Four test cases for the breakage equation are used to acquire closed form series solutions. Further, numerical findings such as number density and moments are compared with the analytical solutions to show the efficiency and accuracy of the method. Moreover, the DJM is employed to solve the well-known ABE, and truncated solutions are presented for the two test cases. In addition, absolute errors over various time periods are depicted in the form of tables.
\keywords{Integro-partial differential equation, Aggregation-breakage equation, DJM, series solution.}
\end{abstract}
\section{Introduction}
Aggregation-breakage equation is model for particle growth caused by aggregation and breakage that transform the dynamics of the particles, i.e., mass or volume. The primary mechanism of these models is the aggregation of two small volume particles to form a bigger particle and the breakage of a particle into tiny particles. These are breakage and aggregation-breakage equations, which are a sort of  integro-partial differential equations.
Integral equations are used to solve a wide range of issues in physics, chemistry, and Biology. As a result, Semi-analytical  techniques to solve integral equations (particularly nonlinear equations) have received more attention in recent years. This paper explains an iterative approach for obtaining solutions to linear and  nonlinear functional equations.

In recent decades, numerical approaches have become the traditional method for solving and evaluating many complicated nonlinear problems. We shall present an alternate strategy in this research, employing iterative technique to obtain a solution with a high degree of accuracy. The iterative approach generates a series that may be summed to discover an analytical formula or used to build an appropriate approximation. By correctly truncating the series, the approximation error may be controlled. \\
The pure breakage equation (\ref{maineq1}) investigates the change in the number density function $c(t,u) \geq 0$ for particles of volume $u> 0$ at time $t \geq 0$ in the physical system as shown in the following equation

\begin{align}\label{maineq1}
\frac{\partial c(t,u)}{\partial t}= \int_{u}^\infty B(u,v)c(t,v)S(v)dv-S(u)c(t,u),
\end{align}
	with the given initial data 
					\begin{align}\label{in1}
						c(0,u)\ \ = \ \ c^{in}(u) \geq 0, \ \ \ u \in ]0,\infty[.
						\end{align}
						
In the Eq.(\ref{maineq1}), the time variable $t$ and the volume variable $u$ are represented as the dimensionless quantities \cite{ramkrishna2000population}.  The selection function $S(u)$ expresses the rate of selecting particles of volume $u$ for breakage into smaller daughter particles.  $B(u, v)$ is the breakage function denoting the formation of daughter particles of volume $u$
from parent particle of volume $v$ due to breakage process. The function $B(u,v)$ is regarded to fulfill the following requirements:
\begin{align}\label{breakagecon1}
    \int_{0}^{v}B(u,v) du=\gamma(v) \geq 1,  \ \ \forall v \geq u,\ \ \ B(u,v)=0, \ \ \forall v < u,
\end{align}
\begin{align}\label{breakagecon2}
    \int_{0}^{v}uB(u,v) du=v. 
\end{align}
In Eq.(\ref{breakagecon1}), $\gamma(v)$ specifies the total amount of fragments formed by a particle of volume $v$. When a parent particle of volume $v$ splits into smaller daughter particles, the total volume of the daughter particles created is equal to $v$. Eq.(\ref{breakagecon2}) demonstrates this assumption known as the conservation of the total volume. \\
We will now  discuss the various terms in Eq.(\ref{maineq1}):
The term on the left side describes the particle number density's temporal rate of change. The first term on the right side of Eq.(\ref{maineq1}) is known as the birth term because it represents the input of particles of volume $u$ to the system as a result of the breakdown of the larger particle of volume $v$. The second one  is death term because it excludes particles of volume $u$ owing to their fragmentation into tiny particles. \\
Now, we have coupled the binary nonlinear aggregation process to the breakage equation (\ref{maineq1}) to get the Eq.(\ref{maineq2}) that is known as the binary aggregation-breakage equation. In binary aggregation, two tiny particles come together to build a large size particle. The mathematical expression of the binary nonlinear  aggregation-breakage equation is:

	\begin{align}\label{maineq2}
						\frac{\partial c(t,u)}{\partial t}\ \ =\ \ &\frac{1}{2}\int_0^{u} K(v,u-v)c(t,v)c(t,u-v)dv-\int_0^\infty K(u,v)c(t,u)c(t,v)dv \nonumber\\
						& +  \int_{u}^\infty B(u,v)S(v)c(t,v)dv-S(u)c(t,u),
					\end{align}
					with the given initial data 
					\begin{align}\label{in2}
						c(0,u)\ \ = \ \ c^{in}(u) \geq 0, \ \ \ u \in ]0,\infty[.
					\end{align}
	In Eq.(\ref{maineq2}), $K(u,v)$ is the aggregation kernel denotes the rate at which particles of volume $u$ and $v$ unite together.  It follows the symmetry condition $K(u,v)=K(v,u)$. The first and second terms on the right side of (\ref{maineq2}) are birth and death terms due to aggregation. These terms illustrate the creation and disappearance of particles of volume  $u.$ The total number of particles fluctuates over time in aggregation-breakage processes while the total mass of particles remains constant. Some integral properties (moments) in terms of $c(t,u)$ have physical interpretation in the particulate system :
	\begin{align}\label{moment}
	    \mu_{j}(t):=\int_{0}^{\infty}{u}^{j}c(t,u) du,\ \ j=0,1,2,....
	    \end{align}
	    $\mu_{0}(t)$ and $\mu_{1}(t)$ describe the total amount of particles and the total volume of particles at time $t\geq 0$.  $\mu_{0}(t)$  is decreased in the aggregation and increased in the breakage processes, while $\mu_{1}(t)$ is a conserved quantity during both processes. $\mu_{2}(t)$  specifies the amount of energy produced by the system. The volume conservation property is the following:
	    \begin{align}\label{massconserve}
	        \int_{0}^{\infty}uc(t,u) du=\int_{0}^{\infty}uc^{in}(u) du, \ \ t\geq 0.
	    \end{align}
	   Several studies have been conducted on the aggregation and breakage equations, and analytical solutions have been discovered for a limited number of examples. Various solutions to the breakage equation have been proposed in \cite{ziff1986kinetics,ziff1991new,ernst1993fragmentation}. 
	    The theoretical results about the existence and uniqueness of solution for Eqs.(\ref{maineq1}) and (\ref{maineq2}) are in \cite{lamb2004existence,laurenccot2000class,mclaughlin1997existence}. There are several numerical methods those are quadrature method of moments \cite{su2007solution,attarakih2009solution}, sectional method \cite{kostoglou2009sectional,kumar2008convergence}, finite element method \cite{ahmed2013stabilized} and finite volume method (FVM) \cite{bourgade2008convergence,kumar2013numerical,forestier2012finite,kumar2015development}.
Bouragde and filbet \cite{bourgade2008convergence} developed the FVM for coupled aggregation and binary breakage equation based on conservative formulation and showed the convergence analysis for locally bounded kernels. The authors in \cite{kumar2013numerical} investigated the FVM for multiple breakage equation and numerically validated the second-order convergence rate.
The authors also discussed the FVM for aggregation equation that conserves the overall volume of the particles in this study \cite{forestier2012finite}. Furthermore, the approach handles overlapping the produced cells by distributing them proportionally throughout the defined meshes. \\
The disadvantages of the numerical schemes are demonstrated by their possible dependency on nonphysical assumptions such as discretization, linearization, sets of basis functions, and many more. Several researchers have recently expressed an interest in semi-analytical techniques to resolve these problems.  Many academics have used semi-analytical techniques to explore the area of aggregation and breakage equations \cite{hammouch2012laplace,singh2015adomian,dutta2018population,kaur2019analytical,hasseine2015two,kaushik2022novel}. Authors \cite{hammouch2012laplace} proposed a modified variational iteration method (VIM) for solving the aggregation equation for constant $(K(u,v)=1)$ and product $(K(u,v)=uv)$ aggregation kernels.  In \cite{singh2015adomian}, the Adomian decomposition method (ADM) is applied to aggregation and breakage equations for several test problems and shows the convergence analysis for particular kernels. For the aggregation and breakage process, the researchers in \cite{kaur2019analytical} developed the homotopy perturbation method (HPM) to obtain the closed form solutions for  $K(u,v)=1, \ u+ v, \ {u}^{2/3}+{u}^{2/3}$ with  $c(o,u)=e^{-u}$ and $K(u,v)=uv$ with $c(o,u)=e^{-u},\frac{e^{-u}}{u}$. In addition, they compared analytical solutions of breakage equation to series solutions for $B(u,v)=2/v, S(v)=v, {v}^2$ with $ c(o,u)=e^{-u},\ \delta(u-r)$ and $B(u,v)=\frac{\alpha+2}{v}\big(\frac{u}{v}\big)^{\alpha},\ S(v)={v}^{\alpha +2}$ with $\ c(o,u)=e^{-u} $ where $-1 < \alpha 
\leq 0.$\\

This work applies the DJM \cite{daftardar2006iterative} over to breakage and aggregation-breakage equations for several test cases. The advantage of DJM is no need to use additional approaches such as estimating Adomian's polynomials as in the ADM or employing Lagrange multipliers in the VIM. Our goal is to discover an approximate solution without generating any restrictive assumptions.  Varsha Daftardar-Gejji and Hossein Jafari were the first to develop the DJM in 2006. This iterative method has effectively solved different types of partial differential equations such as \cite{bhalekar2008new,bhalekar2006solving,yaseen2013exact}. Therefore, we can employ DJM to solve the problems (\ref{maineq1}) \& (\ref{maineq2}) to get the appropriate approximate solution with a high degree of accuracy.  In this article, we have compared the series solution to the breakage equation's analytical solution, and it yields relevant results with lower absolute error. 
Hence, we applied the method over the aggregation-breakage equation (\ref{maineq2}) without having analytical solutions.\\

This article is presented as follows: Section {\ref{iterative method}} discusses the general methodology of DJM. In Section  \ref{NM1}, the plots for moments, particle density function and absolute error are designed for the four different types of kernels in the breakage equation. In addition, error distribution is also shown for various time levels. In Section \ref{NM2}, the ABE is formulated by DJM and demonstrates the particle density function, moments and absolute error plots for two test cases. Finally, conclusions are made in Section \ref{con}.
\section{Iterative method}\label{iterative method}
To show the processes of applying the DJM, we will first analyse the following general functional equation (\ref{eq1}).
\begin{align}\label{eq1}
    c=f+L(c)+N(c),
\end{align}
where $f$ is a given function, $L$ and $N$ are linear and nonlinear operators from a banach space to banach space. 	DJM solution for the Eq.(\ref{eq1}) has the form:
	\begin{align}\label{eq2}
	c=\sum_{i=0}^{\infty}c_{i}.
	\end{align}
	Applying DJM solution on linear and nonlinear operators yield
	\begin{align}\label{eq3}
	    L\left(\sum_{i=0}^{\infty}c_{i}\right)=\sum_{i=0}^{\infty}L(c_{i}).
	\end{align}
		Since $L$ is linear. Nonlinear term decomposed by the following way
		\begin{align}\label{eq4}
		    N\left(\sum_{i=0}^{\infty}c_{i}\right)&=N(c_{0})+\sum_{i=1}^{\infty}\bigg\{N\left(\sum_{j=0}^{i}c_{j}\right)-N\left(\sum_{j=0}^{i-1}c_{j}\right)\bigg\} \nonumber \\
		    &=\sum_{i=0}^{\infty}G_i.
		\end{align}
		Using equations (\ref{eq2}), (\ref{eq3}) and (\ref{eq4}) in Eq.(\ref{eq1}), we get
		\begin{align}\label{sol}
		    \sum_{i=0}^{\infty}c_i=f+\sum_{i=0}^{\infty}L(c_i)+\sum_{i=0}^{\infty}G_i.
		\end{align}
		The DJM series terms are:
		\begin{align}\label{eq5}
		    c_0=f,\,\, \noindent c_{m+1}=L(c_m)+G_{m}, \,\,\noindent m=0,1,2,...   \,.
		\end{align}
		Equation (\ref{sol}) provides the solution for equation (\ref{eq1}) as $ \sum_{i=0}^{\infty}c_i$ , where $c_i$ is derived from equation (\ref{eq5}). An approximate solution ($k$-term) to the problem is represented as follows
	\begin{align}\label{eq6}
			c=\sum_{i=0}^{k-1}c_{i}.
	\end{align}
			
\section{Numerical implementation for breakage}\label{NM1}
In this part, we solve equations (\ref{maineq1}-\ref{in1}) analytically using the DJM. We shall use numerous test problems to illustrate the DJM's efficiency and accuracy. All symbolic and numerical computations are performed with the 'MATHEMATICA SOFTWARE' Package. The DJM recommends that we first rewrite equations (\ref{maineq1}-\ref{in1}) in operator form as follows:
\begin{align}
	L(c(t,u))=\int_{u}^\infty B(u,v)S(v)c(t,v)dv-S(u)c(t,u),		
\end{align}
$$N(c(t,u))=0.$$
\begin{example}
Consider Eqs.(\ref{maineq1}-\ref{in1}) for $S(v)=v$, $B(u,v)=2/v$, and $c(0,u)=e^{-u}.$
\end{example}
Here ${\mathcal{L}}=\frac{\partial}{\partial t }$ and with the assistance of DJM, we get the recursion fromula as:
\begin{align}\label{recursion1}
    c_{0}(t,u)&=e^{-u},\nonumber\\
    c_{m+1}(t,u)&={\mathcal{L}}^{-1}\Big(\int_{u}^\infty \frac{2}{v}c_{m}(t,v) {v}dv-{u}c_{m}(t,u)\Big), \,\,\,\, m=0,1,2,... .
\end{align}
We derive the components $c_m(t, u), m\geq 1$, using the recursive approach (\ref{recursion1}) as follows:
\begin{align*}
    c_{1}(t,u)=&\frac{(-t)^{1}{u}^{-1}e^{-u}}{1!}({u}^{2}-2u+0)\\
     c_{2}(t,u)=&\frac{(-t)^{2}{u}^{0}e^{-u}}{2!}({u}^{2}-4u+2)\\
     \vdots&\\
    c_{m}(t,u)=&\frac{(-t)^{m}{u}^{m-2}e^{-u}}{m!}({u}^{2}-2mu+m(m-1)).
\end{align*}
Now, series solution for this problem is obtained by taking sum of $n$-terms as
\begin{align}\label{problem1sol}
    \Phi_{n}(t,u)=\sum_{m=0}^{n}\frac{(-t)^{m}{u}^{m-2}e^{-u}}{m!}({u}^{2}-2mu+m(m-1)).
\end{align}
Taking the limit over $n$, equation (\ref{problem1sol}) approaches to exact solution of the problem \cite{ziff1985kinetics}
\begin{align*}
    \lim_{n\to \infty}\Phi_{n}(t,u)=(1+t^{2})e^{-u(1+t)}.
\end{align*}

\begin{figure}[htb!]
\centering
\subfigure[DJM Solution]{\includegraphics[width=0.45\textwidth,height=0.35\textwidth]{{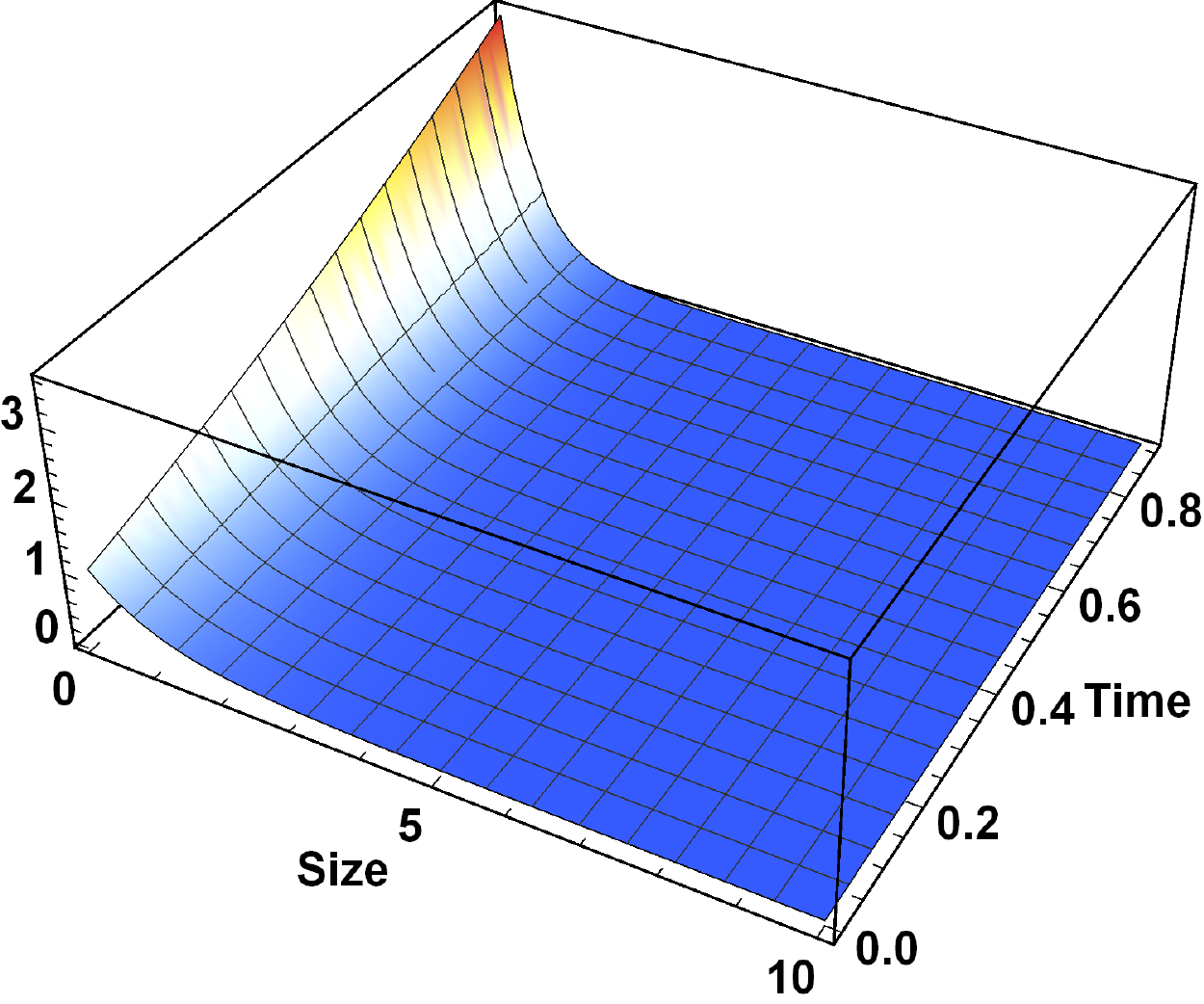}}}
\subfigure[Exact Solution]{\includegraphics[width=0.45\textwidth,height=0.35\textwidth]{{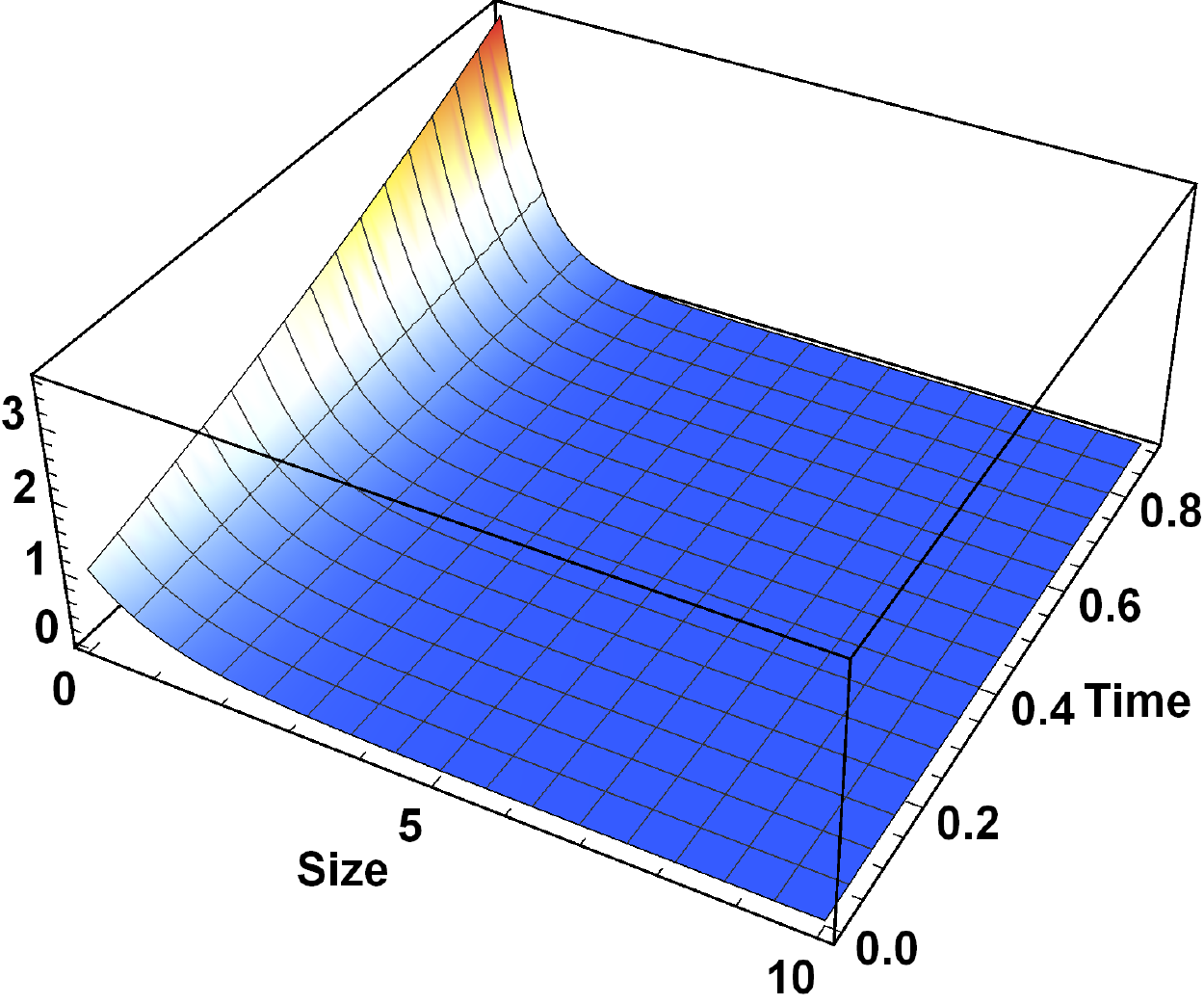}}}
\caption{Approximated and Exact solutions}
\label{fig0}
\end{figure}
\begin{figure}[htb!]
\centering
\subfigure[Absolute Error]{\includegraphics[width=0.45\textwidth,height=0.35\textwidth]{{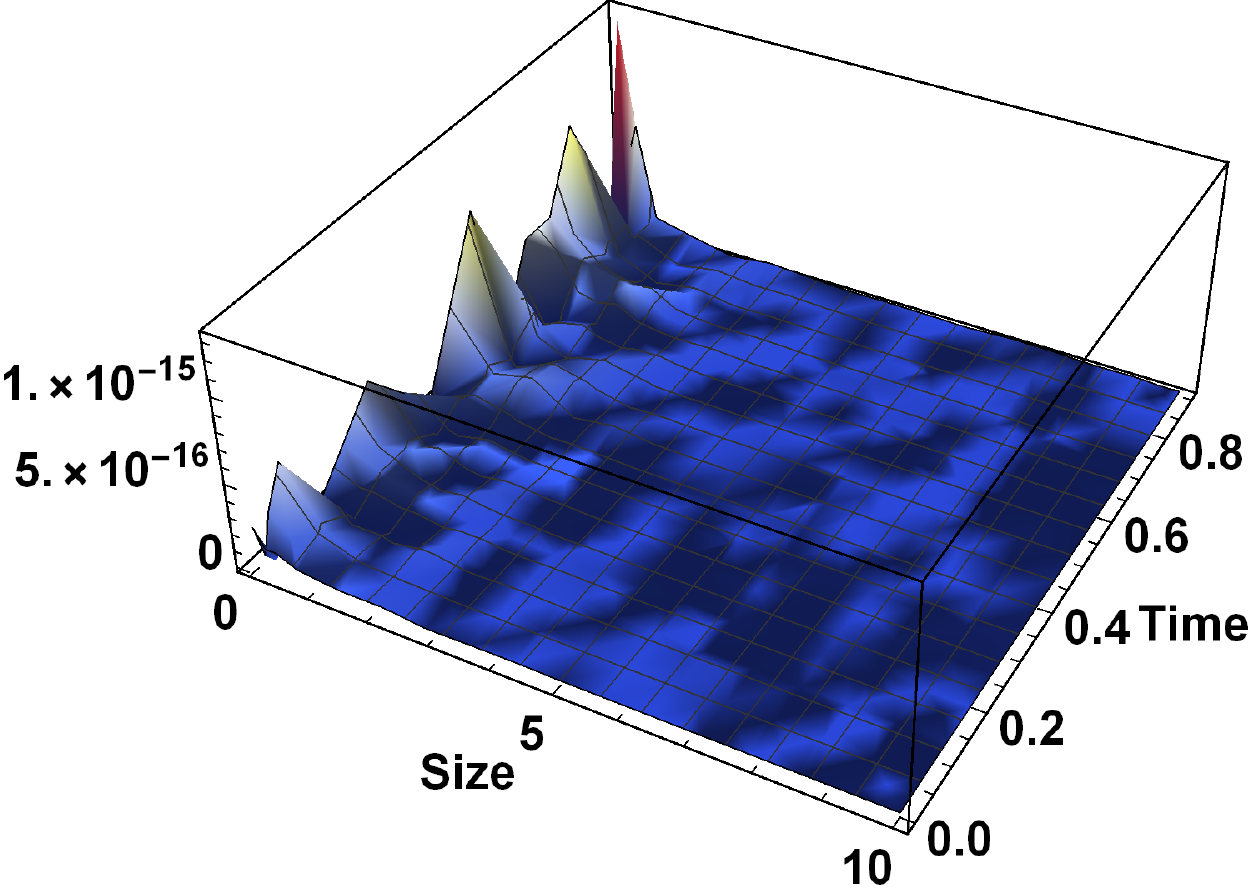}}}
\subfigure[Moments]{\includegraphics[width=0.45\textwidth,height=0.35\textwidth]{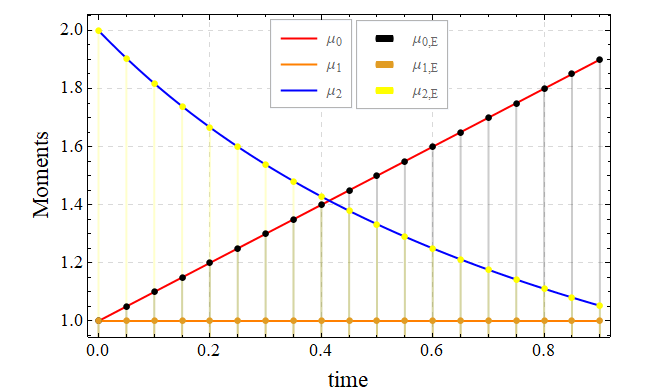}}
 \caption{Absolute error and Moments}
\label{fig1}
\end{figure}
 Figures $1(a)$ and $1(b)$ illustrate the exact solution and 100 terms approximated solution of DJM those are plotted for time 0.9. Considering the DJM solution, one discovers it has a remarkable agreement with the exact solution.  In Figure $2(a)$, the plot shows the maximum absolute error $1. \times 10^{-15}$ that is negligible. Additionally, In Figure $2(b)$, we compared exact moments to approximate solution moments by taking 100 terms, where $\mu_0, \mu_1 $ and $\mu_2$ are the zeroth moment (total number of particles), first moment (total volume of particles) and second moment (energy produced by the system), respectively. $\mu_{0,E},  \mu_{1,E}$ and $ \mu_{2,E}$ are the exact moments. The DJM moments delivered almost equal results and offered an accurate approximation of exact moments.

\begin{table}[h]
    \centering
\begin{tabular}{ |p{0.5cm}||p{2.2cm}|p{2.2cm}|p{2.2cm} |p{2.2cm}|}
 \hline
\multirow{2}{*} {$t$} &\multicolumn{4}{|c|}{Number of terms} \\
 \cline{2-5}
  &10 &15 &20 &25 \\ \hline
0.4 & 1.4536$\times10^{-6}$ & 4.3645$\times10^{-9}$ & 3.5666$\times10^{-12}$ & 9.3987$\times10^{-16}$ \\\hline
0.8 &2.4600$\times10^{-3}$ & 2.3380$\times10^{-4}$& 6.3994$\times10^{-6}$ &4.1016$\times10^{-8}$  \\\hline
1.2&1.8065$\times10^{-1}$ &1.2957$\times10^{-1}$ &2.7863$\times10^{-2}$ &1.3849$\times10^{-3}$  \\\hline
1.6 & 3.7087& 1.1167$\times 10^{-1}$ & 1.0387$\times 10^{-1}$ & 2.2125$\times 10^{-2}$\\\hline
\end{tabular}
\vspace{0.2cm}
 \caption{Error distribution at different level of time for $n$=10,15,20,25.}
    \label{tab:my_label}
\end{table}
Table 1 shows the DJM error at various time levels for n=10,15,20, and 25. As can be seen, the inaccuracy worsens with time for a certain number of terms, whereas the error decreases with more terms evaluated.
\begin{example}
Consider Eqs.(\ref{maineq1}-\ref{in1}) for $S(v)={v}^2$, $B(u,v)=2/v$, and $c(0,u)=e^{-u}.$
    \end{example}
    With the assistance of DJM, we have the following recursion fromula as:
\begin{align}\label{recursion2}
    c_{0}(t,u)&=e^{-u},\nonumber\\
    c_{m+1}(t,u)&={\mathcal{L}}^{-1}\Big(\int_{u}^\infty \frac{2}{v}c_{m}(t,v) {v^2}dv-{u^2}c_{m}(t,u)\Big), \,\,\,\, m=0,1,2,....
\end{align}
Now, having the components $c_m(t, u), m\geq 1$, utilizing the recursive approach (\ref{recursion2}) as follows:
\begin{align*}
    c_{1}(t,u)=&\frac{(-t)^{1}{u}^{0}e^{-u}}{1!}({u}^{2}-2u-2)\\
     c_{2}(t,u)=&\frac{(-t)^{2}{u}^{2}e^{-u}}{2!}({u}^{2}-4u-4)\\
      \vdots&\\
      c_{m}(t,u)=&\frac{(-t)^{m}{u}^{2m-2}e^{-u}}{m!}({u}^{2}-2mu-2m)).
\end{align*}
Series solution for this problem is shown in the following equation
\begin{align}\label{problem2sol}
    \Phi_{n}(t,u)=\sum_{m=0}^{n}\frac{(-t)^{m}{u}^{2m-2}e^{-u}}{m!}({u}^{2}-2mu-2m)).
\end{align}
 Equation (\ref{problem2sol}) yields the exact solution of the problem \cite{ziff1985kinetics} as $n \to \infty$
\begin{align*}
    \lim_{n\to \infty}\Phi_{n}(t,u)=(1+2t+2tu)e^{-u(1+tu)}.
\end{align*}

\begin{figure}[htb!]
\centering
\subfigure[DJM Solution]{\includegraphics[width=0.45\textwidth,height=0.35\textwidth]{{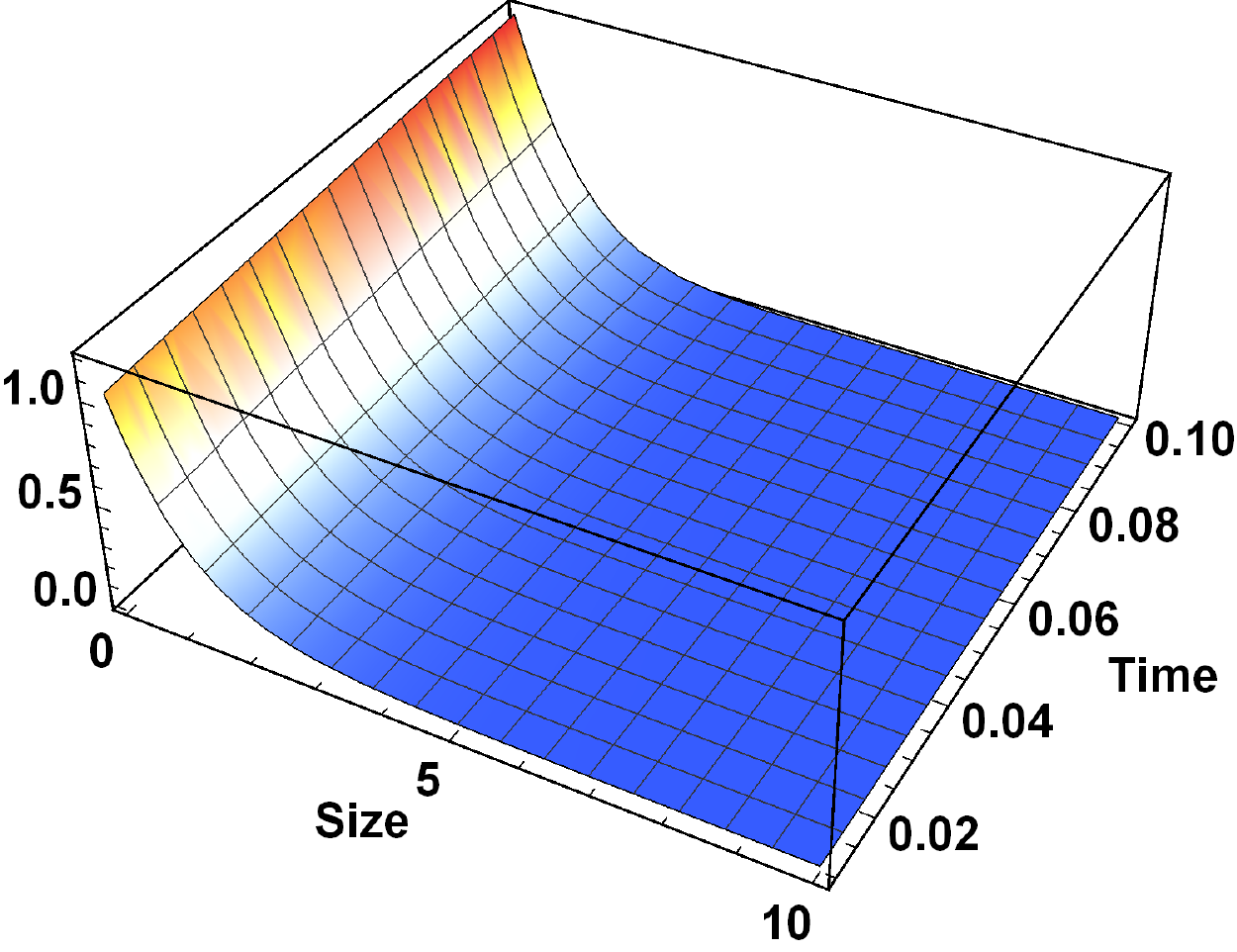}}}
\subfigure[Exact Solution]{\includegraphics[width=0.45\textwidth,height=0.35\textwidth]{{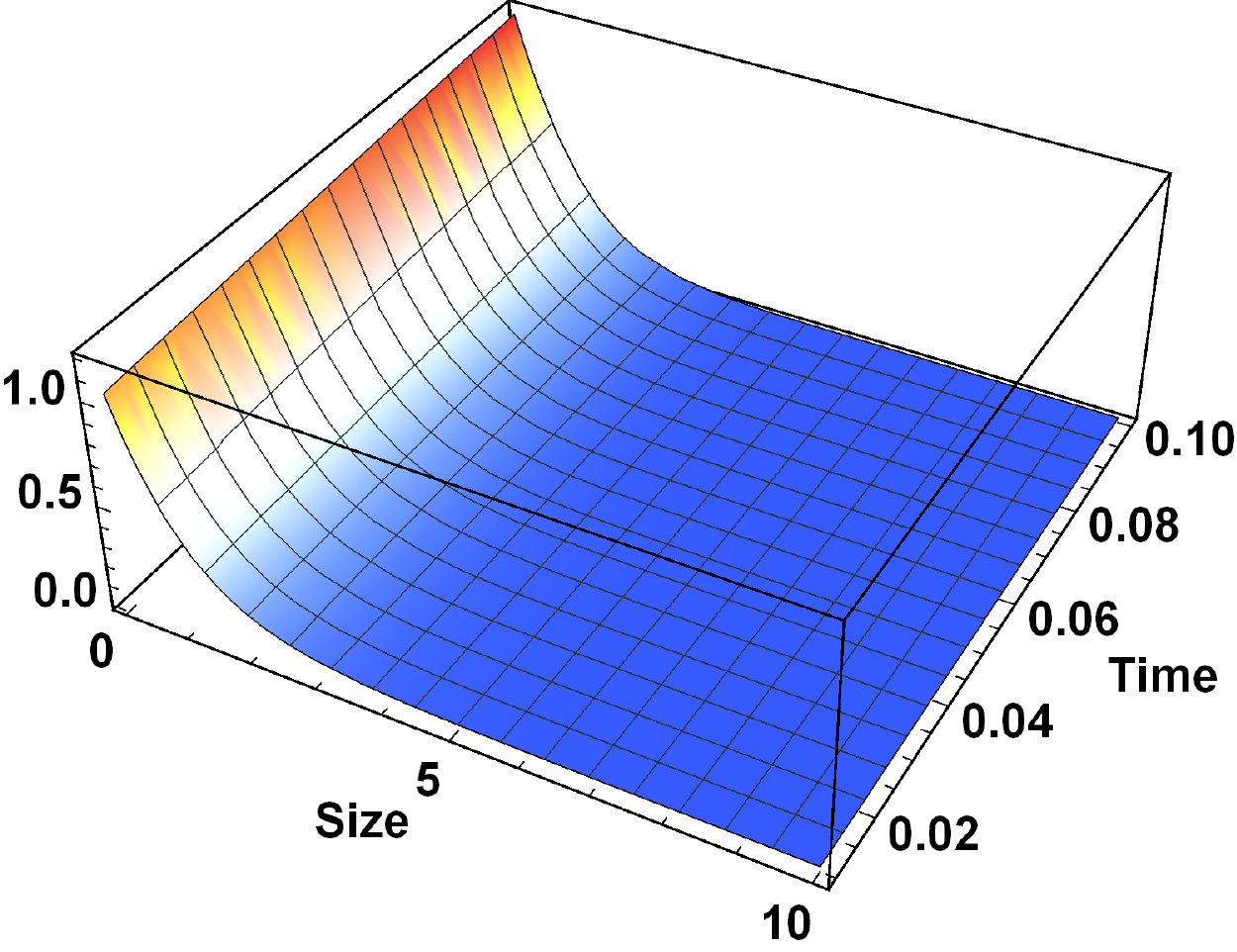}}}
\subfigure[Absolute Error]{\includegraphics[width=0.45\textwidth,height=0.35\textwidth]{{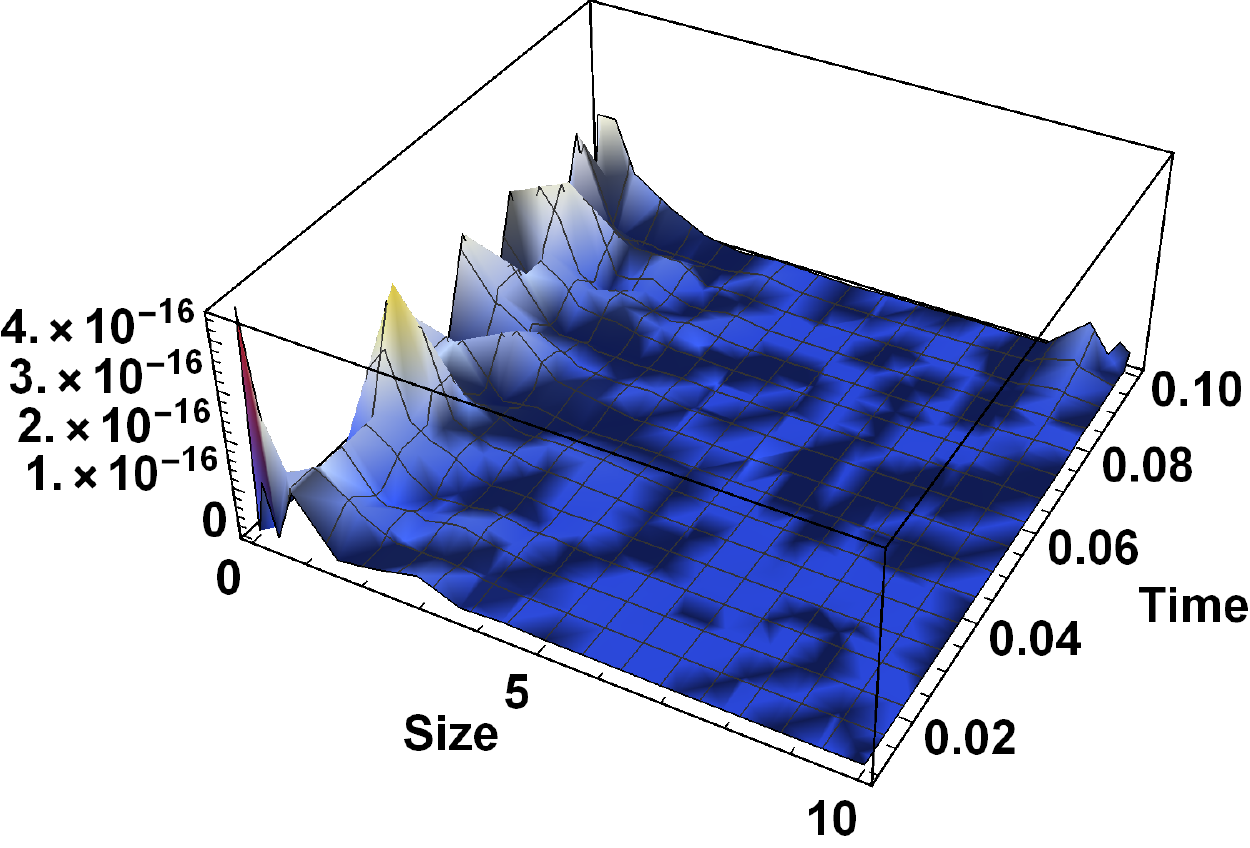}}}
\subfigure[Moments]{\includegraphics[width=0.45\textwidth,height=0.35\textwidth]{{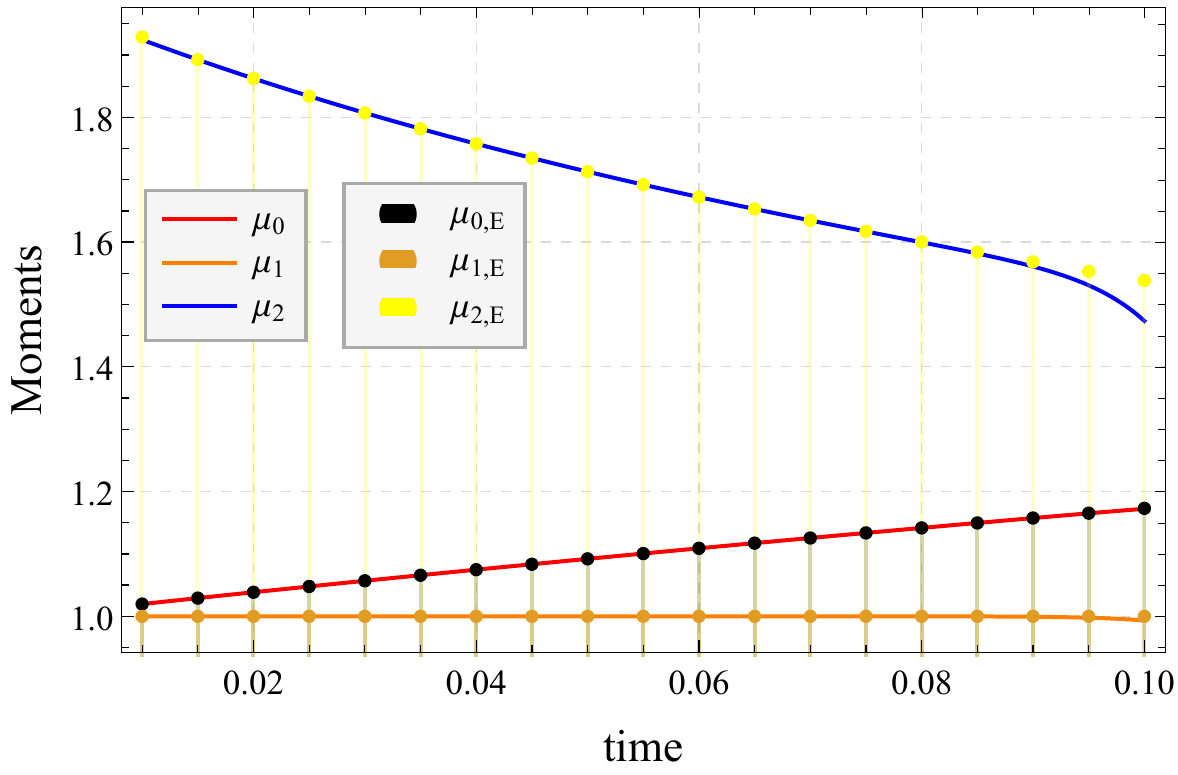}}}
 \caption{Approximated and Exact solutions, Absolute error and Moments}
\label{fig2}
\end{figure}
Figures $3(a)$ and $3(b)$ depict the exact solution and 20 terms approximation solution of DJM for time 0.1. When looking at the DJM solution, one can see that it has a remarkable agreement with the precise solution. The curve in Figure $3(c)$ indicates the highest absolute error $4. \times 10^{-16}$, which is negligible. In addition, we compared exact moments to approximate solution moments in Figure $3(d)$ by choosing 20 terms. The DJM zeroth and first moments produced almost identical findings and provided an accurate estimate of exact moments. However, DJM second moment starts slip away from the exact moment around time 0.1.

 \begin{table}[h]
    \centering
\begin{tabular}{ |p{0.6cm}||p{2.2cm}|p{2.2cm}|p{2.2cm} |p{2.2cm}|}
 \hline
\multirow{2}{*} {$t$} &\multicolumn{4}{|c|}{Number of terms} \\
 \cline{2-5}
  &10 &15 &20 &25 \\ \hline
0.01&1.2396$\times 10^{-12}$ & 1.5304$\times 10^{-16}$ & 1.5082$\times 10^{-16}$ &1.5081$\times 10^{-16}$ \\\hline
0.04 &4.2130$\times 10^{-6}$ & 8.5656$\times 10^{-9}$&3.6889$\times 10^{-12}$ &6.5285$\times 10^{-16}$  \\\hline
0.07&1.6646$\times 10^{-3}$ &5.7854$\times 10^{-5}$ &4.1945$\times 10^{-7}$ &9.2436$\times 10^{-10}$  \\\hline
0.10 & 7.2341$\times 10^{-2}$& 1.5432$\times 10^{-2}$ & 6.7948$\times 10^{-4}$ & 4.5887$\times 10^{-8}$\\\hline
\end{tabular}
\vspace{0.2cm}
 \caption{Error distribution at different level of time for $n$=10,15,20,25.}
    \label{tab:my_label1}
\end{table}
The absolute error is reported in Table \ref{tab:my_label1} for different time levels. As predicted, the error approaches zero as the number of terms in the series solution increases.
\begin{example}
Consider Eqs.(\ref{maineq1}-\ref{in1}) for $S(v)=v$, $B(u,v)=2/v$, and $c(0,u)=\delta({u-r)}.$
\end{example}
According to DJM, we get the recursion fromula as:
\begin{align}\label{recursion3}
    c_{0}(t,u)&=\delta({u-r)},\nonumber\\
    c_{m+1}(t,u)&={\mathcal{L}}^{-1}\Big(\int_{u}^\infty \frac{2}{v}c_{m}(t,v) {v}dv-{u}c_{m}(t,u)\Big), \,\,\,\, m=0,1,2,....
\end{align}
The solution components  $c_m(t, u), m\geq 1$ are determined using the relation (\ref{recursion3})
as
\begin{align*}
    c_{1}(t,u)=&\delta({u-r)}\frac{(-ut)^{1}}{1!}+2t\theta(r-u)\\
     c_{2}(t,u)=&\delta({u-r)}\frac{(-ut)^{2}}{2!}+2t\theta(r-u)\frac{(-ut)^{1}}{1!}+t^2(r-u)\theta(r-u)\\
     c_{3}(t,u)=&\delta({u-r)}\frac{(-ut)^{3}}{3!}+2t\theta(r-u)\frac{(-ut)^{2}}{2!}+\frac{(-ut)^{1}}{1!}t^2(r-u)\theta(r-u)\\
      \vdots &\\
      c_{m}(t,u)=&\delta({u-r)}\frac{(-ut)^{m}}{m!}+2t\theta(r-u)\frac{(-ut)^{m-1}}{(m-1)!}+\frac{(-ut)^{m-2}}{(m-2)!}t^2(r-u)\theta(r-u),
\end{align*}
where $\delta(u-r)$ denotes Dirac's delta function, and $\theta(r-u)$ indicates the unit step function. As a result, the $n$-terms truncated series solution is as follows:
\begin{align}\label{problem3sol}
    \Phi_{n}(t,u)=& \sum_{m=0}^{n}\delta({u-r)}\frac{(-ut)^{m}}{m!}+\sum_{m=0}^{n}2t\theta(r-u)\frac{(-ut)^{m-1}}{(m-1)!} \nonumber \\
    &+\sum_{m=0}^{n}\frac{(-ut)^{m-2}}{(m-2)!}t^2(r-u)\theta(r-u).
\end{align}
 The exact solution of the Example  \cite{ziff1985kinetics} is obtained by, as $n \to \infty$
\begin{align*}
    \lim_{n\to \infty}\Phi_{n}(t,u)=e^{-ut}[\delta(u-r)+\theta(r-u)(2t+t^2(r-u))].
\end{align*}

\begin{example}
Consider Eqs.(\ref{maineq1}-\ref{in1}) for $S(v)={v}^2$, $B(u,v)=2/v$, and $c(0,u)=\delta({u-r)}.$
\end{example}
According to  DJM, we have the following recursive relation:
\begin{align}\label{recursion4}
    c_{0}(t,u)&=\delta({u-r)},\nonumber\\
    c_{m+1}(t,u)&={\mathcal{L}}^{-1}\Big(\int_{u}^\infty \frac{2}{v}c_{m}(t,v) {v}^2dv-{u}^2c_{m}(t,u)\Big), \,\,\,\, m=0,1,2,....
\end{align}
The solution components  $c_m(t, u), m\geq 1$ are determined using the relation (\ref{recursion4})
as
\begin{align*}
    c_{1}(t,u)=&\delta({u-r)}\frac{(-{u}^2t)^{1}}{1!}+2rt\theta(r-u)\\
     c_{2}(t,u)=&\delta({u-r)}\frac{(-{u}^2t)^{2}}{2!}+2rt\theta(r-u)\frac{(-{u}^2t)^{1}}{1!}\\
     c_{3}(t,u)=&\delta({u-r)}\frac{(-{u}^2t)^{3}}{3!}+2rt\theta(r-u)\frac{(-{u}^2t)^{2}}{2!}\\
      \vdots & \\
      c_{m}(t,u)=&\delta({u-r)}\frac{(-{u}^2t)^{m}}{m!}+2rt\theta(r-u)\frac{(-{u}^2t)^{m-1}}{(m-1)!}.
\end{align*}
Here, $ \Phi_{n}(t,u)$ expresses the $n$-terms truncated series solution of the problem:
\begin{align}\label{problem4sol}
    \Phi_{n}(t,u)= \sum_{m=0}^{n}\delta({u-r)}\frac{(-{u}^2t)^{m}}{m!}+\sum_{m=0}^{n}2rt\theta(r-u)\frac{(-{u}^2t)^{m-1}}{(m-1)!}.
\end{align}
We obtained the exact solution of the Example \cite{ziff1985kinetics} as $n \to \infty$
\begin{align*}
    \lim_{n\to \infty}\Phi_{n}(t,u)=e^{-{u}^2t}[\delta(u-r)+2rt\theta(r-u)].
\end{align*}
\section{Numerical results for ABE}\label{NM2}
All these results related to the breakage equation (\ref{maineq1}) have been concluded with four examples, and we have the exact solutions for the examples. DJM has given efficient results. Consequently, it would be a suitable method for the ABE (\ref{maineq2}) to obtain the approximate solution. In this article, we will examine two test cases of ABE that lack analytic solutions.
Firstly, Eqs.(\ref{maineq2}-\ref{in2}) can be transformed into an operator expression after applying DJM as
\begin{align}
	L(c(t,u))=\int_{u}^\infty B(u,v)S(v)c(t,v)dv-S(u)c(t,u),		
\end{align}
\begin{align}
N(c(t,u))= \frac{1}{2}\int_0^{u} K(v,u-v)c(t,v)c(t,u-v)dv-\int_0^\infty K(u,v)c(t,u)c(t,v)dv.\end{align}

\begin{example}
Consider Eqs.(\ref{maineq2}-\ref{in2}) for $K(x,y)=1$, $S(v)=v$, $B(u,v)=2/v$  and $c(0,u)=e^{-u}.$
\end{example}
Applying the DJM approach, the following recursive relation is obtained
\begin{align}\label{recursion5}
    c_{0}(t,u)=& e^{-u},\nonumber\\
    c_{m+1}(t,u)=&{\mathcal{L}}^{-1}\Big(\int_{u}^\infty \frac{2}{v}c_{m}(t,v) {v}dv-{u}c_{m}(t,u)\Big) \nonumber \\ 
    &+{\mathcal{L}}^{-1}\Big(\frac{1}{2}\int_{0}^{u} \sum_{k=0}^{m}c_{m}(t,u-v) \sum_{k=0}^{m}c_{m}(t,v)dv \nonumber \\
    &-\int_{0}^{\infty}\sum_{k=0}^{m}c_{m}(t,u) \sum_{k=0}^{m}c_{m}(t,v)dv\Big) \nonumber \\
    &+{\mathcal{L}}^{-1}\Big(\frac{1}{2}\int_{0}^{u} \sum_{k=0}^{m-1}c_{m}(t,u-v) \sum_{k=0}^{m-1}c_{m}(t,v)dv\nonumber\\
    & -\int_{0}^{\infty}\sum_{k=0}^{m-1}c_{m}(t,u) \sum_{k=0}^{m-1}c_{m}(t,v)dv\Big)
   \quad  m=0,1,2,....
\end{align}
The solution components $c_m(t,u), m \geq 1$ are determined using the relation (\ref{recursion5}) and due to the complexities of the terminology, just a few are written here as
\begin{align*}
    c_{1}(t,u)=& t \ (2e^{-u}-e^{-u}u)+t \ \Big(-e^{-u}+\frac{e^{-u}u}{2}\Big)\\
     c_{2}(t,u)=&\frac{1}{4}e^{-u}t^2(2-4u+{u}^2)-t \ \Big(-e^{-u}+\frac{e^{-u}u}{2}\Big)\\&+\frac{1}{144}e^{-u}t(72(-2+u)-18t(6+(-6+u)u)+t^2(-24+(-6+u)^2 u)).
\end{align*}
A five term truncated solution is considered using the recursive technique defined in equation (\ref{recursion5}). In Figure $4(a)$, the DJM offered an approximated solution of 5 terms at time 0.4. It is noted in figure $4(b)$ that the approximated moments ($\mu_{j,5}, j=0,1,2$) with 5 terms and approximated moments ($\mu_{j,4}, j=0,1,2$) with 4 terms are almost  identical. 
\begin{figure}[htb!]
\centering
\subfigure[Density function  (n=5)]{\includegraphics[width=0.45\textwidth,height=0.35\textwidth]{{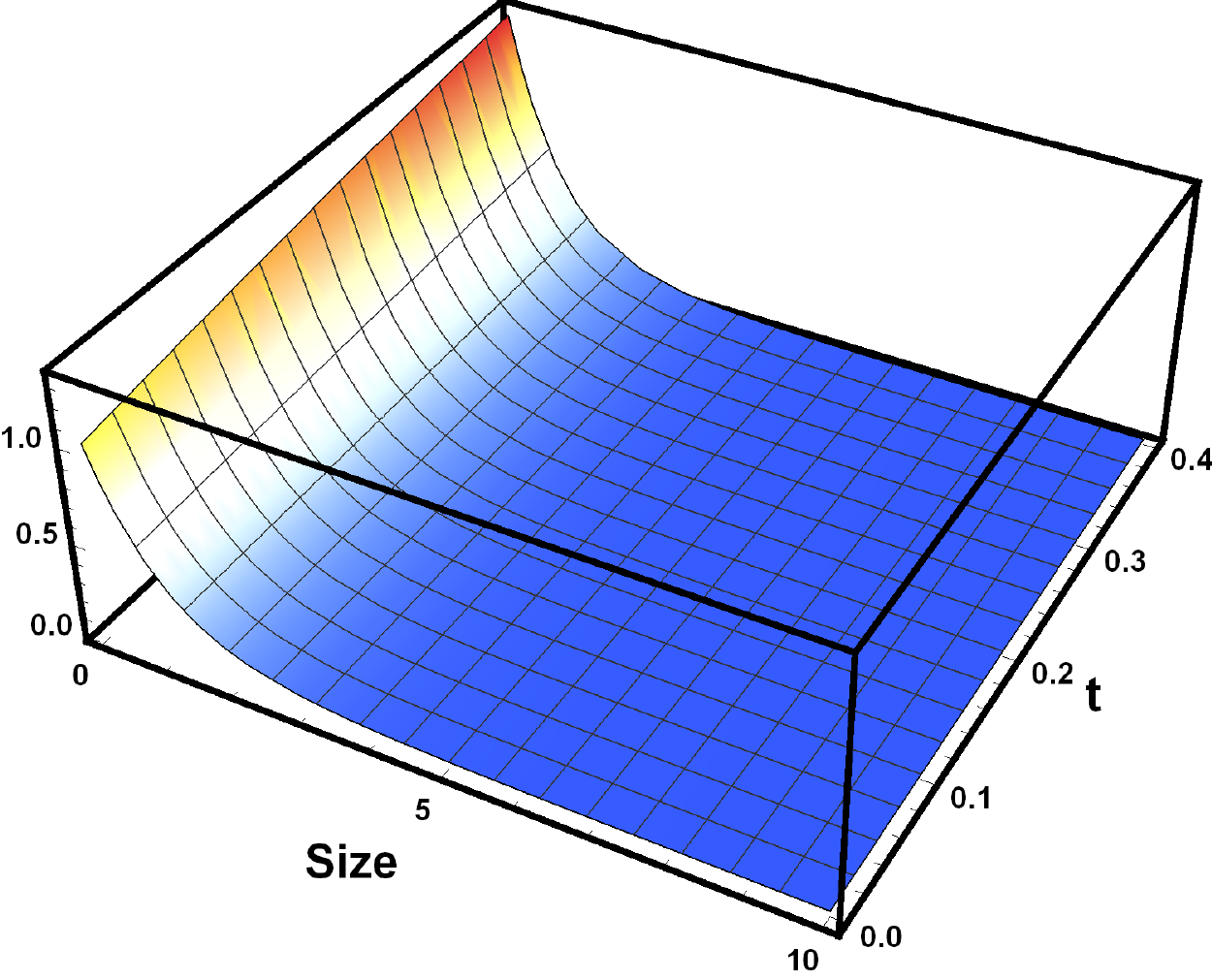}}}\label{density5}
\subfigure[Moments]{\includegraphics[width=0.45\textwidth,height=0.35\textwidth]{{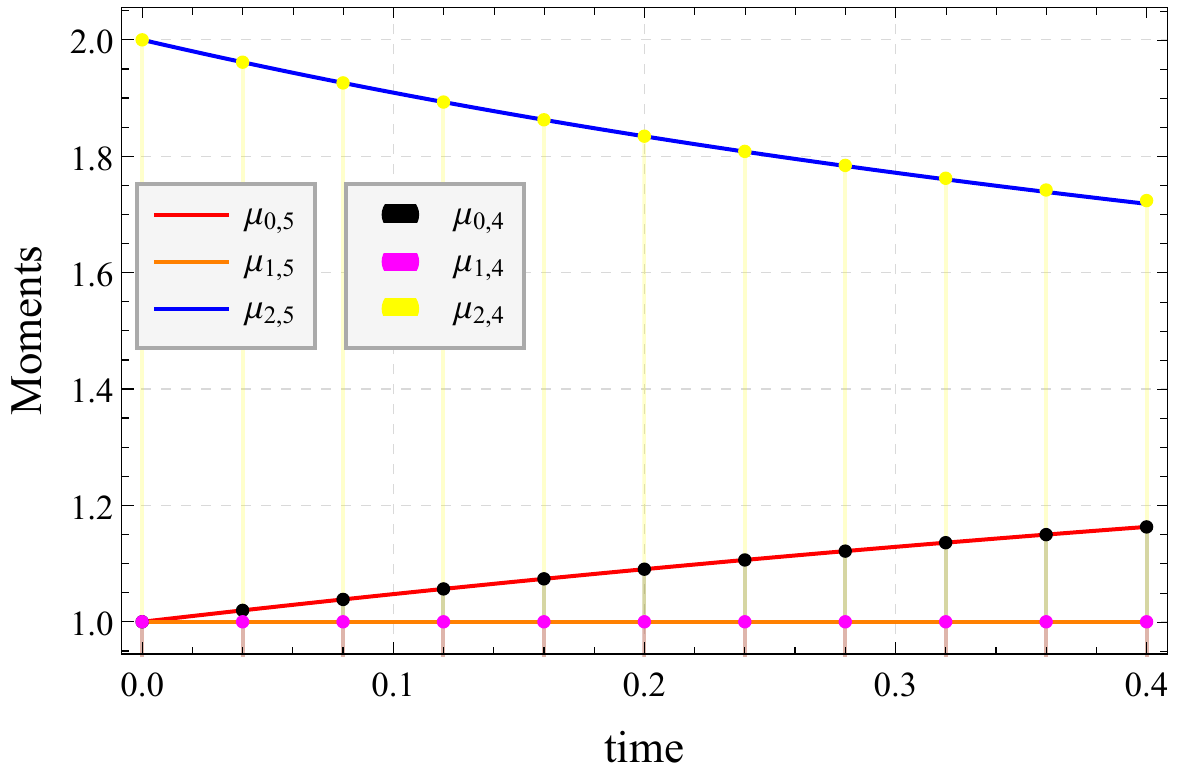}}}\label{moments5}
 \caption{Density function and Momentss}
\label{fig3}
\end{figure}
Figure $5(a),5(b),5(c)$ and $5(d)$ demonstrate the absolute error behavior among the 1, 2, 3, 4 terms and 5 terms of the series solution as we can observe that the error is going to decrease for a particular time.
\begin{figure}[htb!]
\centering
\subfigure[Absolute error 1]{\includegraphics[width=0.45\textwidth,height=0.35\textwidth]{{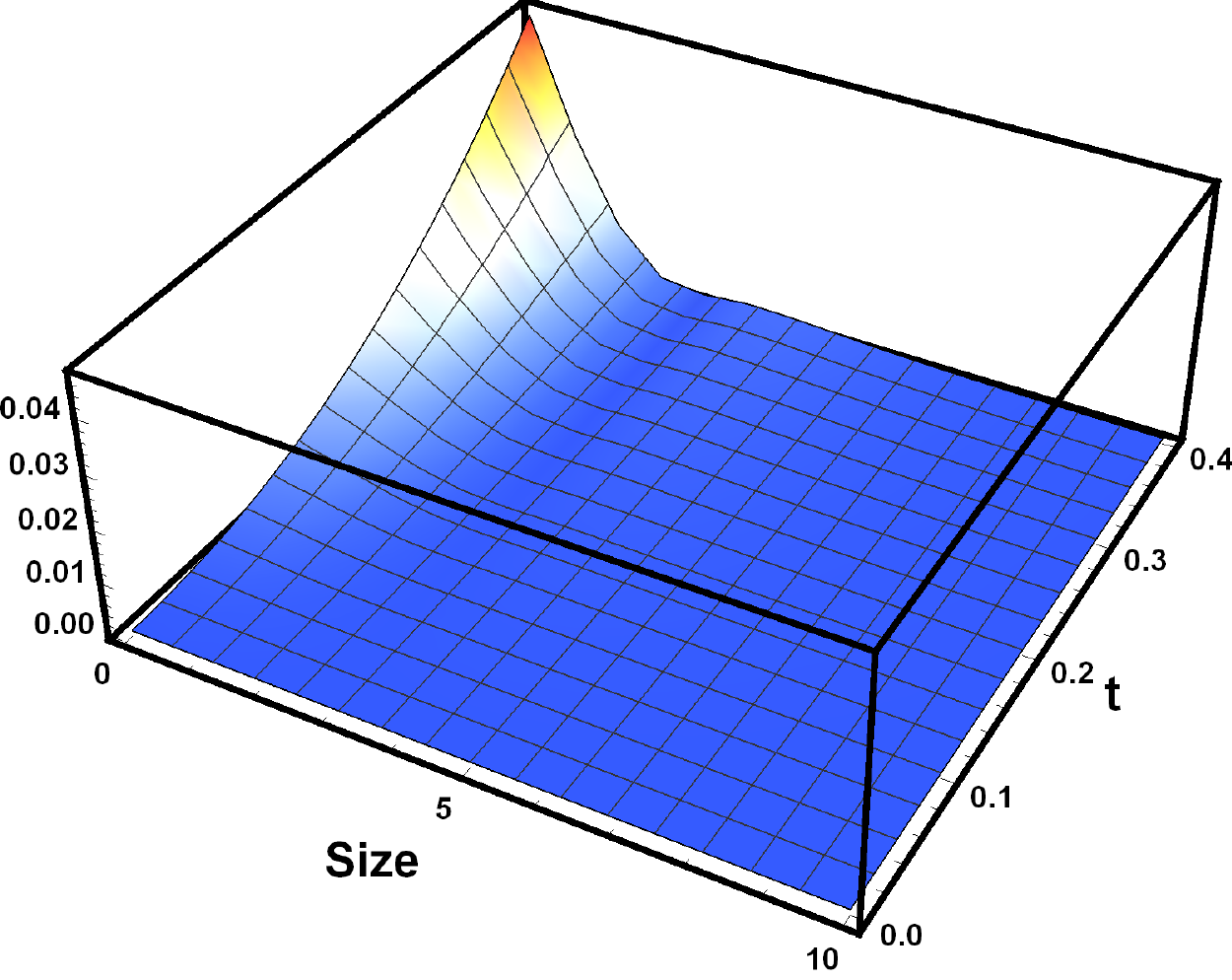}}}
\subfigure[Absolute error 2]{\includegraphics[width=0.45\textwidth,height=0.35\textwidth]{{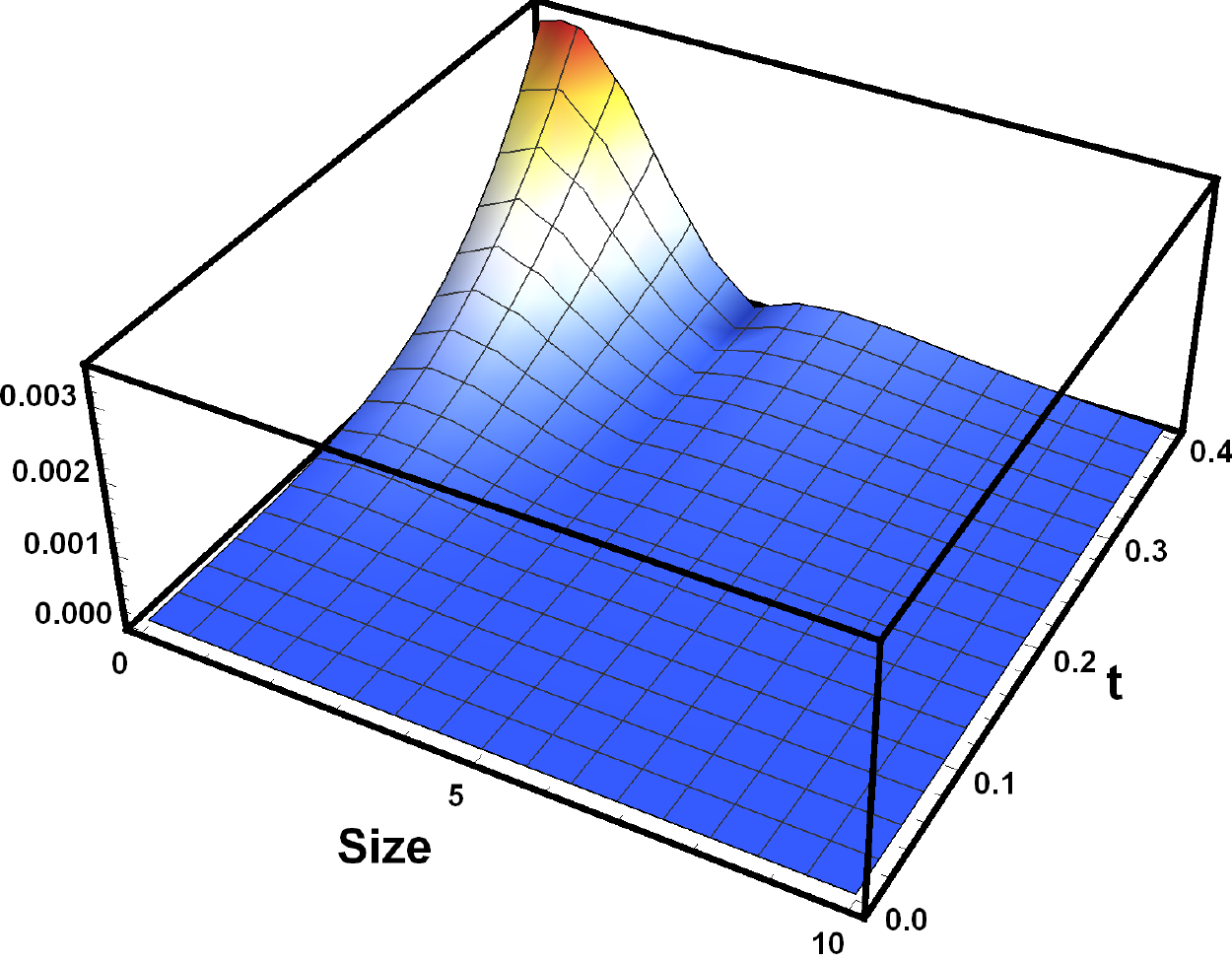}}}
\subfigure[Absolute error 3]{\includegraphics[width=0.45\textwidth,height=0.35\textwidth]{{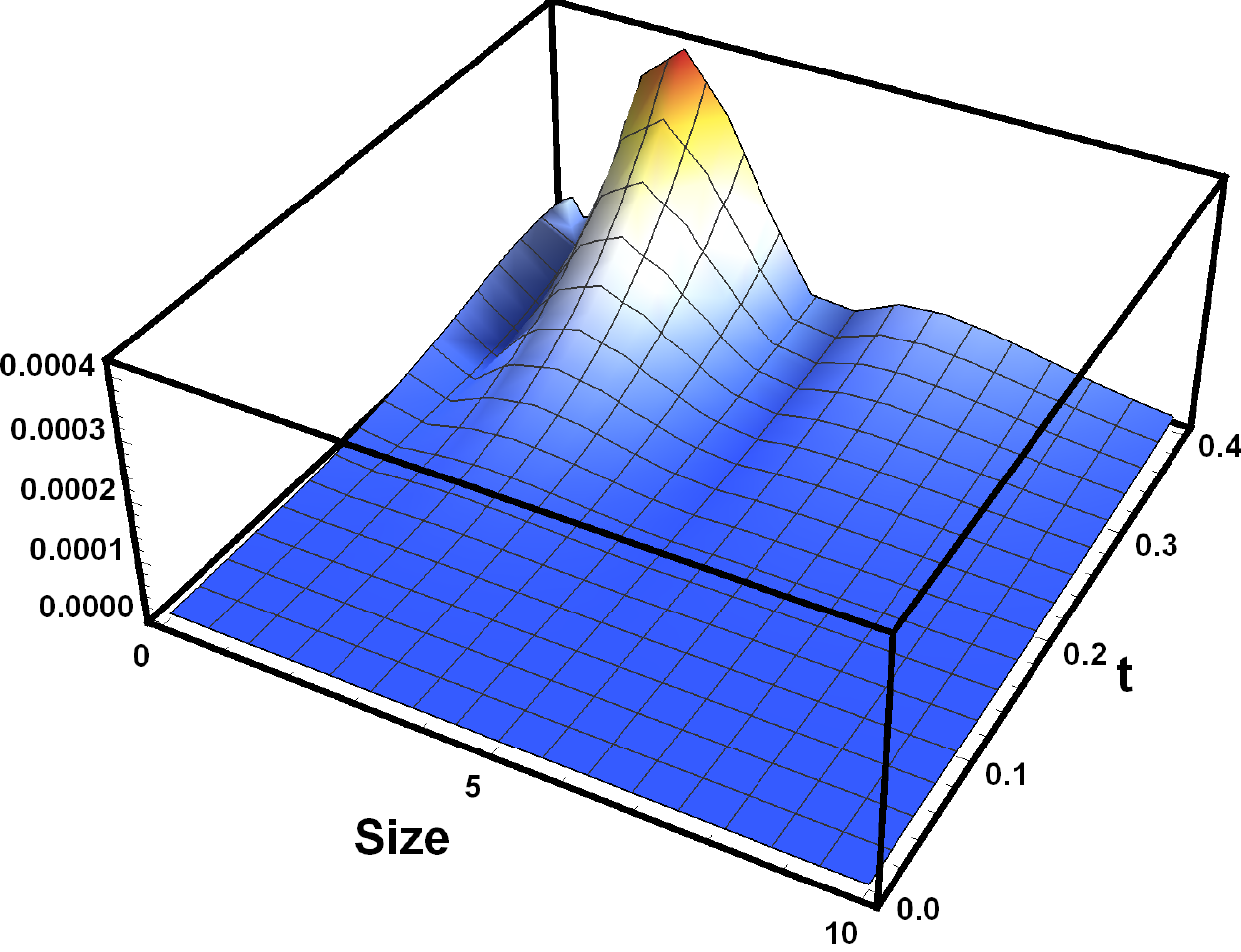}}}
\subfigure[Absolute error 4]{\includegraphics[width=0.45\textwidth,height=0.35\textwidth]{{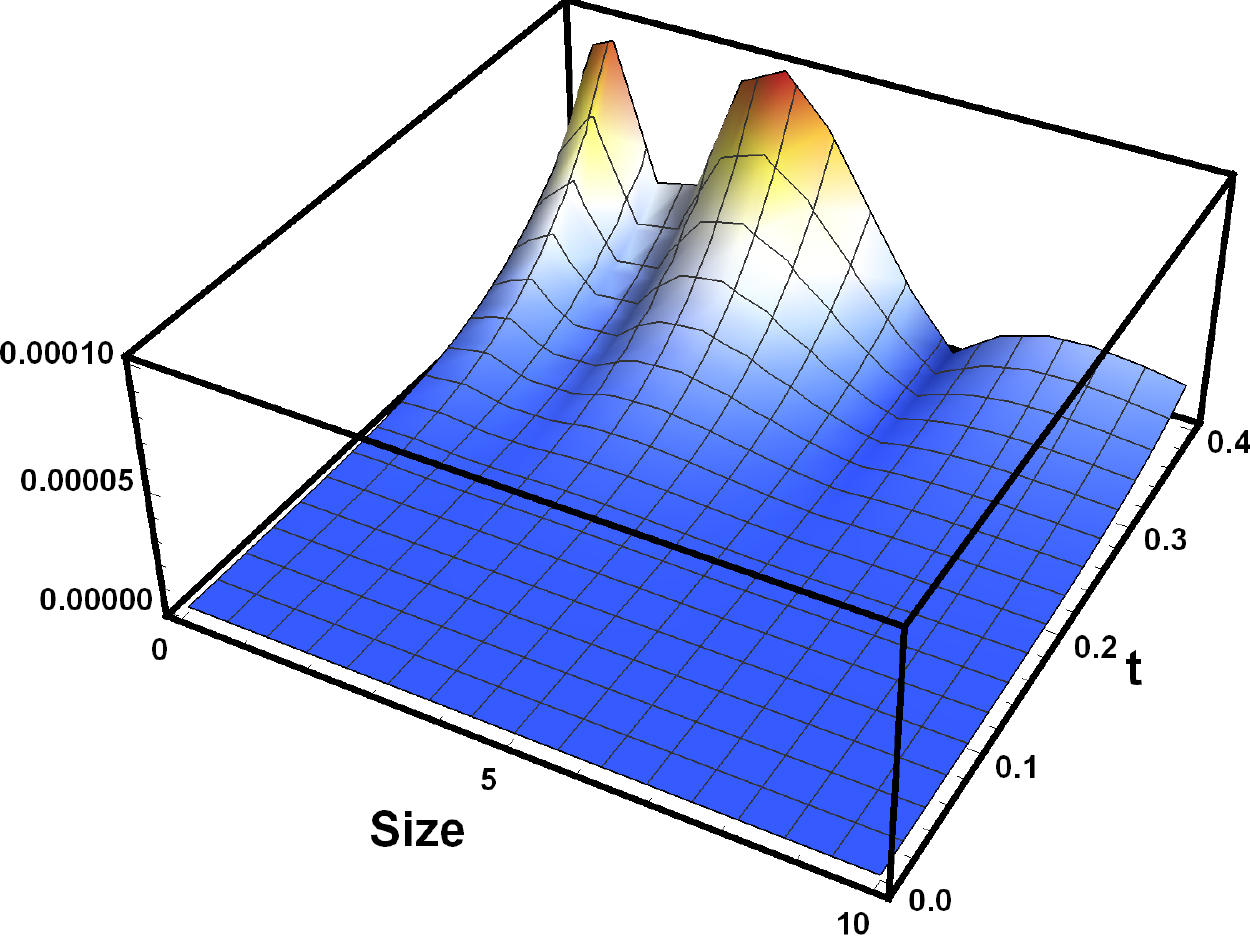}}}
 \caption{ Absolute errors}
\label{fig4}
\end{figure}

\begin{example}
Consider Eqs.(\ref{maineq2}-\ref{in2}) for $K(x,y)=1$, $S(v)=v$, $B(u,v)=2/v$  and $c(0,u)=4ue^{-2u}.$
\end{example}
We got the same recursive relation as in (\ref{recursion5}) except first term $c_0(t,u)=4ue^{-2u}$. The approximate solution components are:
\begin{align*}
    c_{0}(t,u)=& 4ue^{-2u}\\
     c_{1}(t,u)=& t\big(-4e^{-2u}u+\frac{4}{3}e^{-2u}{u}^3\big)+t\big(-4e^{-2u}{u}^2+e^{-2u}(2+4u)\big)\\
      c_{2}(t,u)=&  \frac{1}{2}e^{-2u}t^2-2e^{-2u}t^2u-e^{-2u}t^2{u}^2+\frac{8}{3}e^{-2u}t^2{u}^3-\frac{2}{3}e^{-2u}t^2{u}^4\\&-t\big(-4e^{-2u}u+\frac{4}{3}e^{-2u}{u}^3\big)+\frac{1}{945}e^{-2u}t(1260u(-3+{u}^2)\\&+63t(-15+u(-15+2u(30+x(-5+(-5+u)u)))))\\&+ t^2(-315+2u(315+u(315+u(-525+u(105+u(42+(-14+u)u)))))).
\end{align*}
The recursive approach in equation (\ref{recursion5}) is used to consider a truncated four terms solution except initial term. The DJM provided an estimated solution of four terms at time 0.4 in Figure $6(a)$. Figure $6(b)$ shows that the approximated moments ($\mu_{j,4}, j=0,1,2$) with four terms and approximated moments ($\mu_{j,3}, j=0,1,2$) with three terms seem to be almost similar.
\begin{figure}[htb!]
\centering
\subfigure[Density function]{\includegraphics[width=0.45\textwidth,height=0.35\textwidth]{{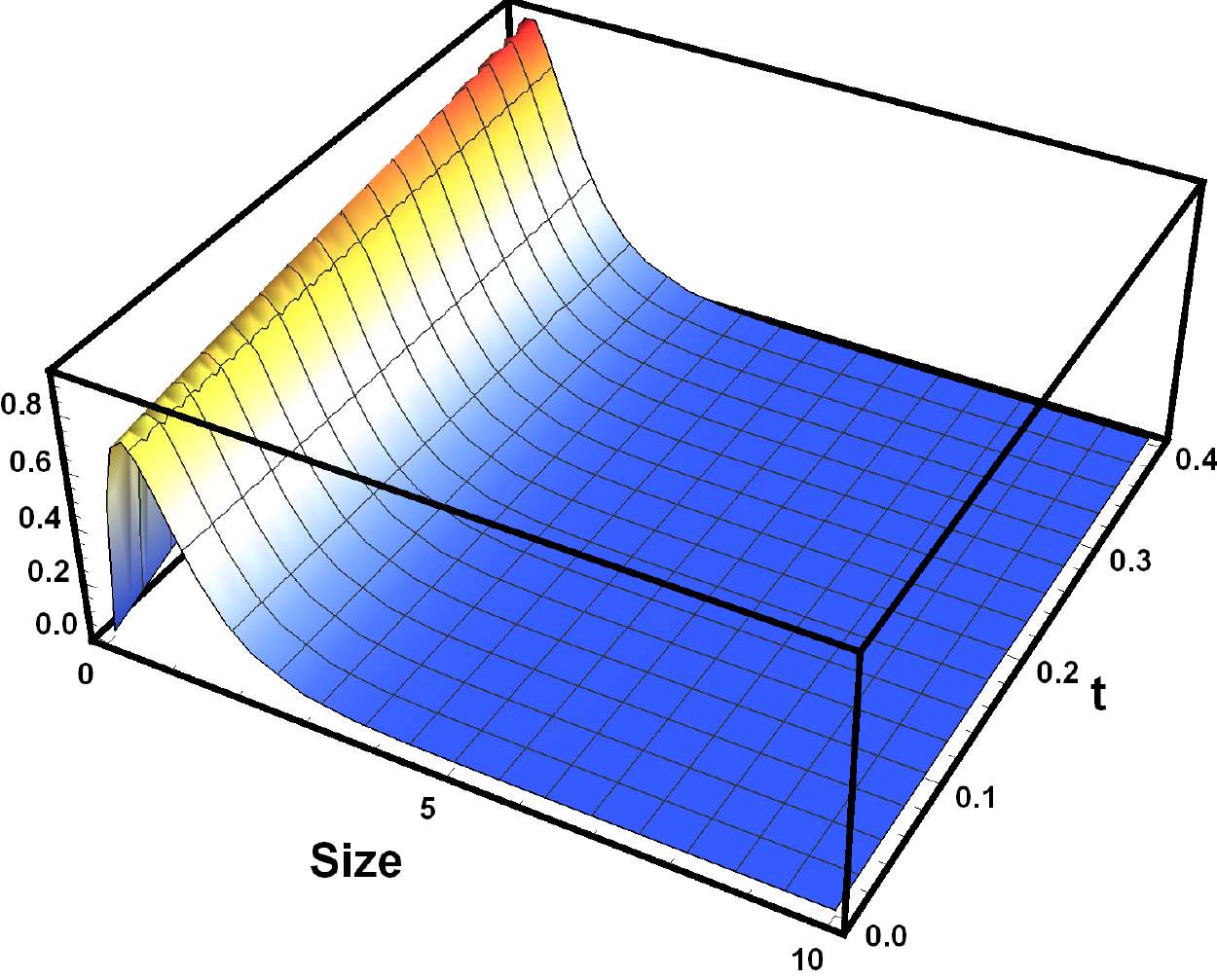}}}
\subfigure[Moments]{\includegraphics[width=0.45\textwidth,height=0.35\textwidth]{{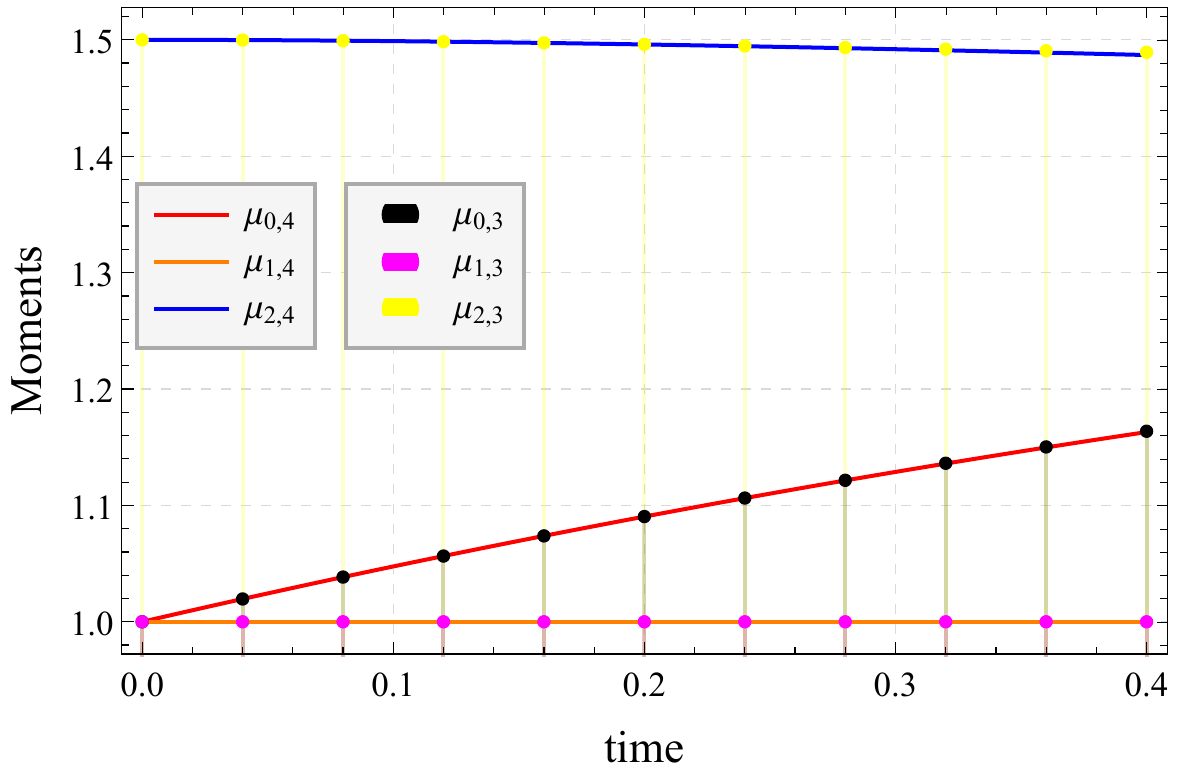}}}
 \caption{Density function and Momentss}
\label{fig5}
\end{figure}

\begin{figure}[htb!]
\centering
\subfigure[Absolute error 1]{\includegraphics[width=0.45\textwidth,height=0.35\textwidth]{{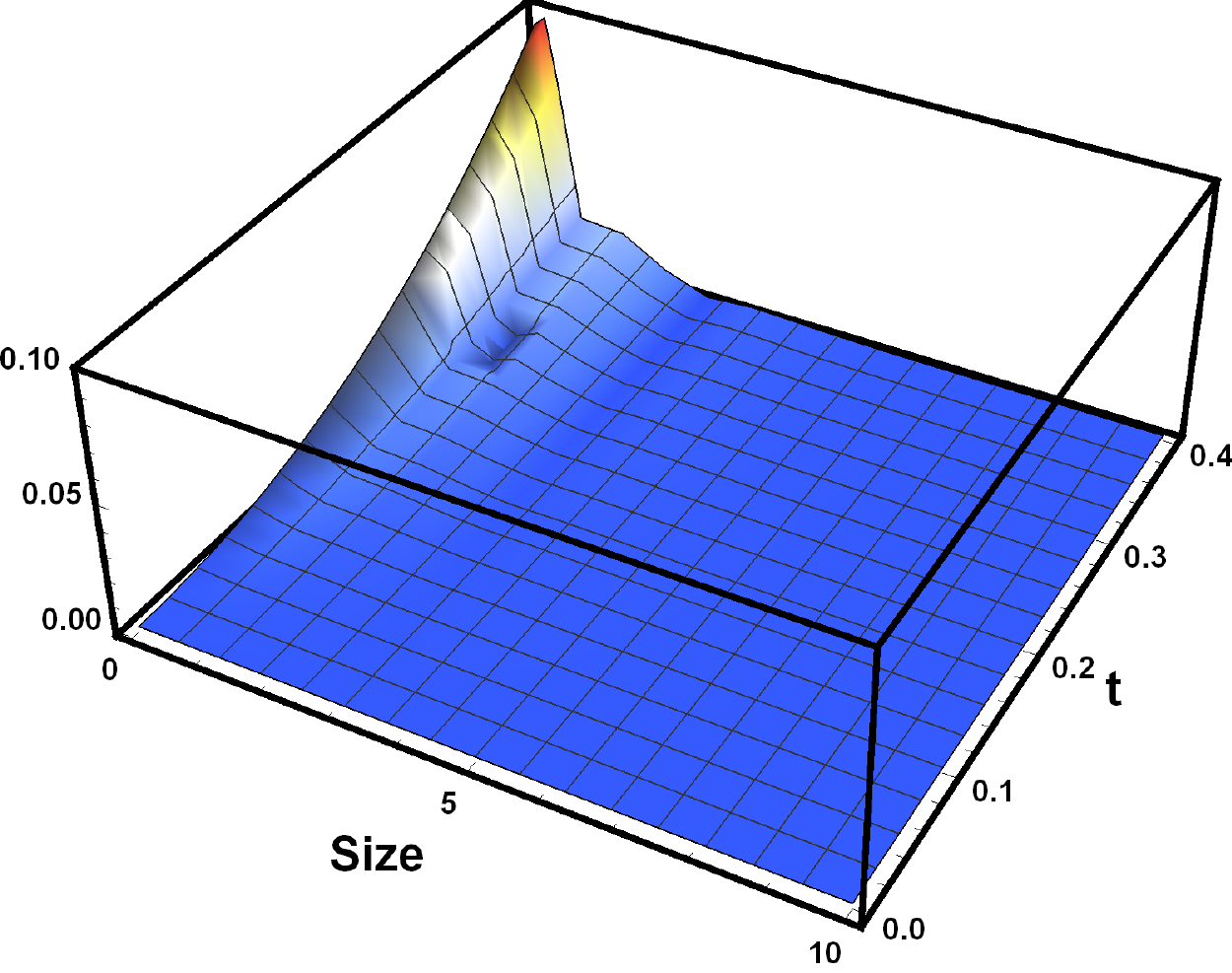}}}
\subfigure[Absolute error 2]{\includegraphics[width=0.45\textwidth,height=0.35\textwidth]{{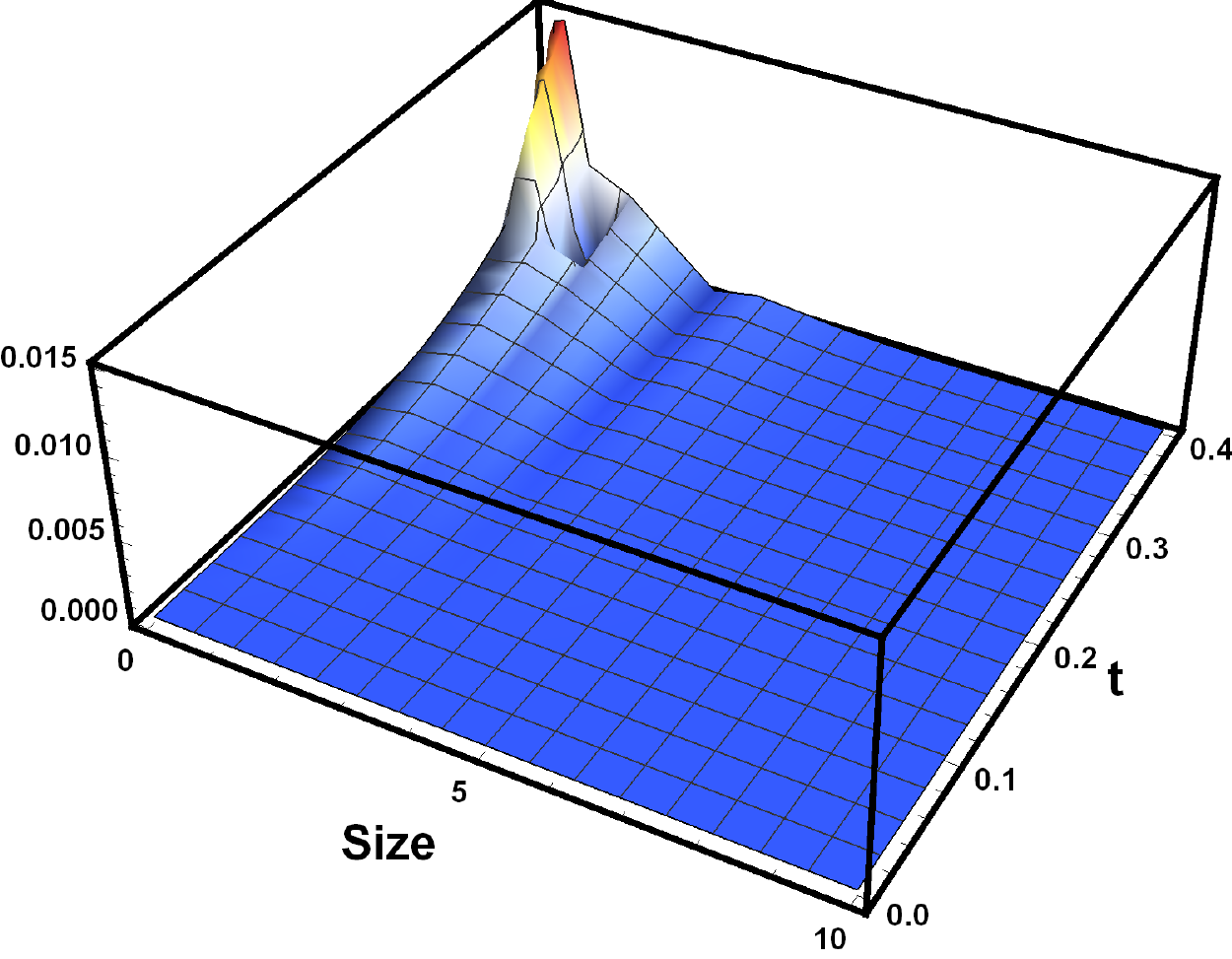}}}
\subfigure[Absolute error 3]{\includegraphics[width=0.45\textwidth,height=0.35\textwidth]{{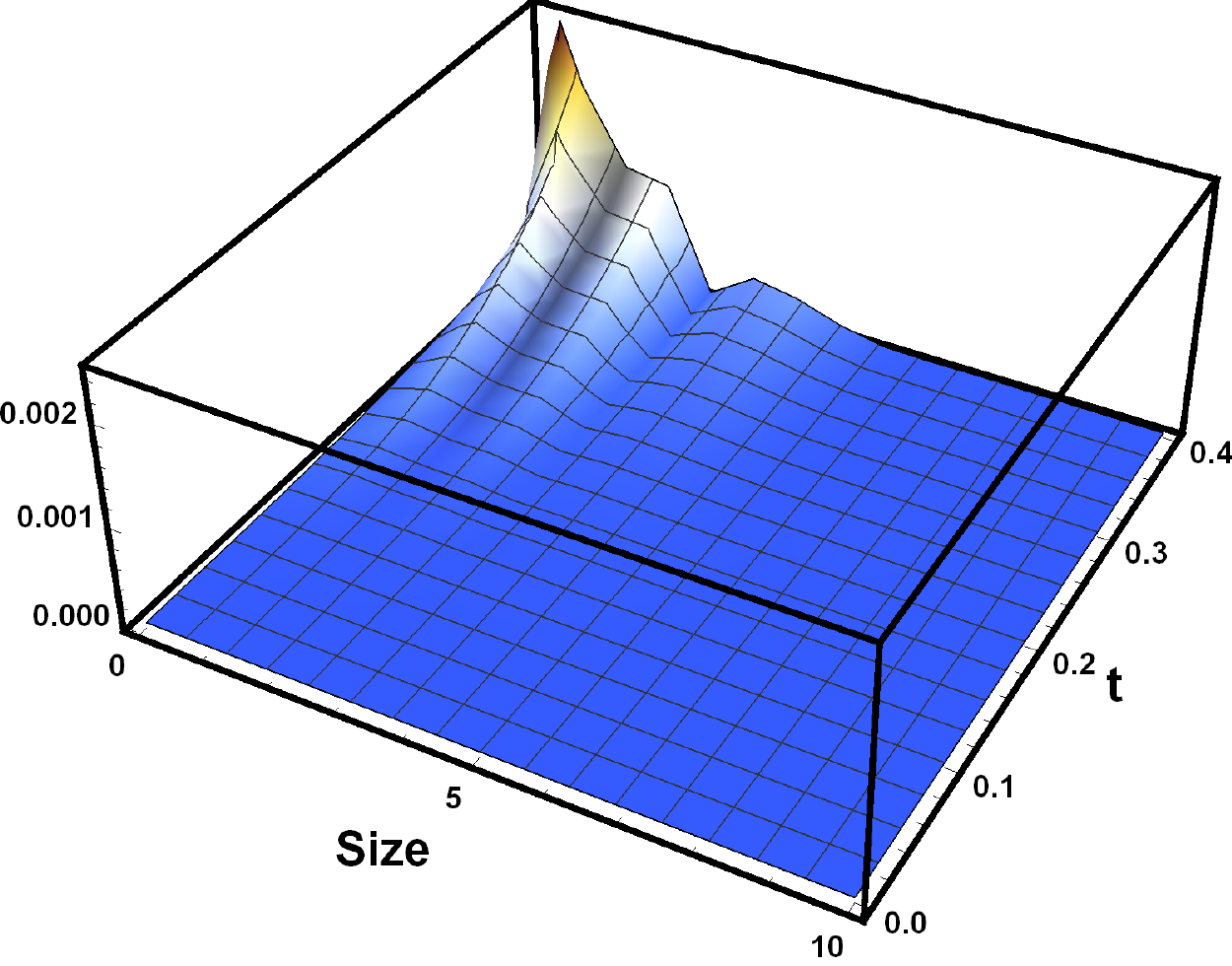}}}
 \caption{ Absolute errors}
\label{fig6}
\end{figure}
In figure $7(a)$, the absolute error between 1 and 4 terms, in figure $7(b)$, the absolute error between 2 and 4 terms and in figure $7(c)$, the absolute error between 3 and 4 terms are shown. The absolute errors among the series solutions are in the decreasing mode at a fixed time.
\section{Conclusions}\label{con}  The purpose of this article was to introduce DJM as a method for solving breakage and aggregation-breakage equations. Four test cases for the breakage equation were evaluated, and closed-form solutions were obtained in every instance. In addition, the accuracy and efficacy of the proposed scheme were demonstrated by comparing truncated and exact solutions, and absolute errors were also presented in tables and graphs. Moreover, moments obtained via DJM exhibited a high degree of precision with regard to the exact moments.  Due to the accuracy acquired in the breakage equation, results are extended for the nonlinear coupled aggregation-breakage equation. In the absence of analytical solutions, we obtained finite-term truncated solutions and plotted number density and moment comparisons.
\clearpage

\bibliographystyle{spmpsci}
\bibliography{DJM}

\end{document}